\documentstyle{article}
\begin{document}
\Large
\begin{center}
{\bf CONICS, ($q+1$)-ARCS, PENCIL CONCEPT OF TIME AND PSYCHOPATHOLOGY}
\end{center}
\vspace*{-0.1cm}
\begin{center}
Metod Saniga\\

\vspace*{.2cm}
\normalsize
Astronomical Institute of the  Slovak Academy of Sciences,\\
05960 Tatransk\' a Lomnica, Slovak Republic\\
(E-mail: msaniga@astro.sk)
\end{center}

\vspace*{.3cm}
\normalsize
\noindent
{\bf  Abstract:} --  It  is  demonstrated  that in the (projective plane over)
Galois  fields  {\sf GF}($q$)  with  $q = 2^{n}$  and  $n \geq 3$ ($n$ being a
positive  integer)  we  can  define,  in  addition  to the temporal dimensions
generated by pencils of conics, also time coordinates represented by aggregates
of  ($q+1$)-arcs  that  are  {\it not}  conics.  The  case is illustrated by a
(self-dual)  pencil  of conics endowed with two  singular  conics of which one
represents a double real line and the other is a real line pair. Although this
pencil  does not generate the ordinary (i.e., featuring  the past, present and
future)  arrow  of  time  over  {\sf GF}($2^{n}$),  there {\it does}  exist  a
pencil-related family of ($q+1$)-arcs, not conics, that closely resembles such
an  arrow.  Some  psycho(patho)logical  justifications  of  this  finding  are
presented,  based  on the `peculiar/anomalous' experiences of time by a couple
of schizophrenic patients.

\vspace*{.1cm}
\\
\\
\large
{\bf 1. Introduction}

\vspace*{.1cm}
\noindent
In  one  of  our  debut  papers  devoted  to  the  theory  of pencil-generated
temporal dimensions [1] we  discussed  basic  properties  of  the structure of
time  over a Galois field of even order, {\sf GF}(2$^{n}$).  Our attention was
focused  exclusively  on  a  specific  pencil of conics featuring two singular
conics and two distinct
base points. Although  this  pencil has been found to reproduce quite well the
qualitative  properties  of the physical world when considered over the ground
field  of the real numbers [2], it leads to a very peculiar arrow of time over
{\sf GF}($2^{n}$), the one lacking (the moment of) the present [1]. The aim of
this short  contribution  is  to  show that there exists an interesting way of
`recovering/restoring' the familiar arrow of time also in the latter case.
\\ \\
{\bf 2. Conics and ($q+1$)-Arcs}

\vspace*{.1cm}
\noindent
To this end in view, let us consider, following (the symbols and notation of)
[1,2], a conic, i.e., the curve of second order
\begin{equation}
{\cal Q}_{\breve{x} \breve{x}} \equiv \sum_{i \leq j} c_{ij} \breve{x}_{i}
\breve{x}_{j} = 0,~~~~~~i, j = 1, 2, 3;
\end{equation}
here  $c_{ij}$  are  regarded  as fixed quantities, while $\breve{x}_{i}$ are
viewed  as  variables  (the so-called homogeneous coordinates of a projective
plane).  The  conic  is  degenerate (or singular) if there exists a change of
coordinate system reducing Eq. (1) into a form in fewer variables; otherwise,
the  conic is non-degenerate (or proper).  It  is well known (see, e.g., [3])
that the equation of any proper conic  ${\cal Q}$  in a projective plane over
{\sf GF}($q$)  (the  latter  being  henceforth  denoted as PG$(2, q)$) can be
brought into the canonical form
\begin{equation}
\widetilde{\cal Q}_{\breve{x} \breve{x}} = \breve{x}_{1} \breve{x}_{2} -
\breve{x}_{3}^{2} = 0.
\end{equation}
From the last equation it follows that the points of $\widetilde{\cal Q}$ can
be parametrized as $\varrho \breve{x}_{i} = (\sigma^{2},1, \sigma)$, $\varrho
\neq 0$, and  this  implies  that a proper conic in PG$(2, q)$ contains $q+1$
points; the point (1, 0, 0) and $q$ other points specified  by  the sequences
$(\sigma^{2}, 1, \sigma)$  as  the  parameter  $\sigma$  runs through the $q$
elements of {\sf GF}($q$).  Moreover,  it  can  easily  be verified  that any
triple of distinct points of $\widetilde{\cal Q}$ are linearly independent,
for [4]
\begin{equation}
\det \left( \begin{array}{ccc}
1 & 0 & 0 \\
\sigma_{1}^{2} & 1 & \sigma_{1} \\
\sigma_{2}^{2} & 1 & \sigma_{2}
\end{array} \right)
= \sigma_{2} - \sigma_{1} \neq 0
\end{equation}
and
\begin{equation}
\det \left( \begin{array}{ccc}
\sigma_{1}^{2} & 1 & \sigma_{1} \\
\sigma_{2}^{2} & 1 & \sigma_{2} \\
\sigma_{3}^{2} & 1 & \sigma_{3}
\end{array} \right)
= (\sigma_{1} - \sigma_{2})(\sigma_{2} - \sigma_{3})
(\sigma_{3} - \sigma_{1}) \neq 0.
\end{equation}
Hence, a proper conic of PG$(2, q)$ is an example of a $(q+1)$-$arc$ -- a set
of $q+1$ points no three of which are collinear. Although every non-degenerate
conic of  PG$(2, q)$ is a ($q+1$)-arc, the converse is true only for $q$ odd;
for $q$ even  and  greater  than  four there also exist ($q+1$)-arcs that are
{\it not} conics [3--5].

In order  to  see this explicitly, we first recall that all the tangents to a
proper conic ${\cal Q}$ of PG$(2, q = 2^{n})$  are concurrent, i.e. they pass
via {\it  one  and  the  same} point, called the nucleus [1,3--5]. Hence, the
conic ${\cal Q}$  together  with  its nucleus form a ($q+2$)-arc. Now, let us
delete  from  this  ($q+2$)-arc  a  point  belonging  to  ${\cal Q}$; we will
obtain  a  ($q+1$)-arc ${\cal K}$ which shares $q = 2^{n}$ points with ${\cal
Q}$. Taking into account  the  fact that a proper conic is uniquely specified
by  five  of its points, no three collinear, it then follows that  the  above
described  ($q+1$)-arc  ${\cal K}$  can$not$  be a conic for $n \geq 3$; for,
indeed,  if  it were then it would have with ${\cal Q}$ more than five points
in common and would thus coincide with the latter, a contradiction.
\\ \\
{\bf 3. Pencil-Time Comprising ($q+1$)-Arcs}

\vspace*{.1cm}
\noindent
By the  very definition, a straight-line (henceforth only line) of PG$(2, q)$
can  have  at  most two points in common with a ($q+1$)-arc ${\cal K}$; if it
has just  two, it is called -- following the terminology used for conics -- a
secant  of  ${\cal K}$,  if one, a tangent to ${\cal K}$, and if none, a line
external to (or, skew with) ${\cal K}$.  So, a ($q+1$)-arc can be regarded as
a  natural  and  straightforward  generalization  of the concept of conic for
PG$(2, 2^{n})$, $n \geq 3$.  As  a  consequence, instead  of viewing the time
dimension as being
generated  by  a  pencil  of conics, we can introduce its generalized concept
based on a one-parametric family of ($q+1$)-arcs.  Moreover, after affinizing
PG$(2, q)$ we define, in a completely analogous way to what we did in the case
of pencils of conics [2], the domain of the past/future to  be represented by
those   ($q+1$)-arcs   (of  a  given  family)  to  which  the ideal  line  is
secant/external,  while  the  ($q+1$)-arc(s)  having  this  line as a tangent
stand(s) for the moment(s) of the present.

In  order to see an explicit realization of this idea, we will again consider
our favoured pencil of conics [1,2]
\begin{equation}
{\cal Q}_{\breve{x} \breve{x}}^{\vartheta} = \vartheta_{1}\breve{x}_{1}
\breve{x}_{2}
+ \vartheta_{2} \breve{x}_{3}^{2} = 0.
\end{equation}
The pencil features two distinct base points, $B_{1}$: $\varrho \breve{x}_{i}
=  (0, 1, 0)$  and  $B_{2}$:  $\varrho  \breve{x}_{i} = (1, 0, 0)$,  each  of
multiplicity  two,  and  two  degenerate  conics:  $\vartheta  \left(  \equiv
\vartheta_{2}/ \vartheta_{1}) \right) = \pm \infty$ (i.e.,  the  double  real
line $\breve{x}_{3}^{2}=0$)  and  $\vartheta = 0$ (i.e., a pair of real lines
$\breve{x}_{1} = 0$  and  $\breve{x}_{2} = 0$  concurring  at  the point $N$:
$\varrho \breve{x}_{i} = (0, 0, 1)$). As already mentioned, this pencil, when
affinized in such a  way that the ideal line $\cal{L}^{\infty}$ meets neither
$B_{1, 2}$ nor $N$, reproduces nicely the ordinary arrow of time if considered
over the field of the reals [2],  but leads to a very peculiar arrow, the one
lacking  the  present, when we switch to {\sf GF}(2$^{n}$) [1]; this  happens
because  the  point $N$ is the common nucleus for {\it all} the proper conics
of pencil (5)  and  as ${\cal L}^{\infty}$ is not incident with $N$ it cannot
be a tangent  to any of them. Let us select one line, ${\cal L}^{\ast}$, from
the pencil of  lines  carried  by $N$ such that ${\cal L}^{\ast} \neq NB_{1},
NB_{2}$. It is  obvious that  the  point $A$ at which ${\cal L}^{\infty}$ and
${\cal L}^{\ast}$  meet each  other  lies  on  just one (proper) conic ${\cal
Q}^{\ast}$  of  pencil  (5),  to  which ${\cal L}^{\infty}$ must clearly be a
secant. Now, let  us  create  a  family of ($q+1$)-arcs in such a way that we
delete  from  each  proper conic the point at which ${\cal L}^{\ast}$ touches
the latter, and add  to  such  a  $q$-arc  the  nucleus  $N$  (recalling once
again that $N$ is the nucleus  for  {\it all}  the proper conics of (5)). The
family of ($q+1$)-arcs  created this  way thus contains not only ($q+1$)-arcs
to which the  ideal line ${\cal L}^{\infty}$ is a secant  (the  past)  and/or
an  external  line  (the  future),  as  in the case of conics [1], but also a
unique ($q+1$)-arc,  that  composed of ${\cal Q}^{\ast} \backslash \{A\}$ and
the point $N$, having with  ${\cal L}^{\infty}$ just one point in common (and
standing thus for the present): this aggregate is thus qualitatively identical
with a geometrical structure that was in [2] recognized and  demonstrated  to
reproduce remarkably well our {\it ordinary/normal} perception of time.
\\ \\
{\bf 4. Some Intriguing Psychopathology of Time}

\vspace*{.1cm}
\noindent
We  have  thus arrived at a principally new type of temporal arrow that cannot
be  reproduced  by  (any pencil of) conics whatever field we would take as the
ground  field  of  the projective plane. And because this kind of the temporal
emerges  only  over  fields  of  characteristic  two  which, we conjecture(d),
represent the `working  regime'  of  the  deepest  parts  of  our subconscious
[1],  the  corresponding  mental  states will be extremely difficult (and, so,
very rare)  to  attain  and  be fully experienced. Nevertheless, after looking
carefully   through  a large  number  of  references  dealing  with  so-called
`altered' states of consciousness [see the exhaustive bibliography in Ref. 6],
we succeeded in finding a very interesting old paper [7] that seems to contain
descriptions of such mental states by schizophrenic patients.
\footnote{The  article  is  written in German; the English translation of both
the excerpts quoted was borrowed from [8] (first excerpt) and [9] (second one). The italics, however, are supplied by the present author.}
Below  are  the  excerpts  from the narratives given by a couple of psychotics
where  there is a/n direct/explicit reference to a `strange,' or `new,' mental
temporal   dimension;   in   particular, a patient, aliased  `Sche,' describes
their `weird' time fabric as follows [7, p. 556]:

\normalsize
\vspace*{.2cm}
\noindent
 ``\ldots  and  then  came  a  feeling of horrible expectation that I could be
sucked  up  into the past or that the past would overcome me and flow over me.
It  was  disquieting  that  someone  could  play with time like that, somewhat
daemonic.  This  would  be  perverse  for humanity. What could time be for the
orderlies?  Did they still  have ordinary time? Then {\it everything seemed to
be  absolutely  of  no consequence}, and in spite of that I was very uneasy. A
{\it foreign}  time  sprang up. Everything was confused, pell-mell, and I felt
contracted  in  myself:  I  wanted  to  hold everything back, but I had to let
everything  go\ldots ~I  wanted  this  {\it  false}  time  to disappear  in me
again\ldots"

\vspace{.2cm}

\large
\noindent
Another  patient  (`Ge')  gives  even a more `physically attractive' piece of
information [7, p. 567]:

\normalsize
\vspace*{.2cm}
\noindent
 ``One  evening  during  a  walk  in a busy street, I had a sudden feeling of
nausea\ldots ~Afterwards  a  small  patch  appeared before my eyes\ldots ~The
patch  glimmered  inwardly  and  there was a to and fro of dark threads\ldots
~The  web  grew more pronounced\ldots ~I felt drawn into it. It was really an
interplay  of movements which had replaced my own person. Time had failed and
stood  still --  no,  it  was  rather  that time {\it re}-appeared just as it
disappeared.  This  {\it new} time was {\it infinitely manifold and intricate
and could hardly be compared with what we ordinarily call time}. Suddenly the
idea  shot through my head that time lies {\it not only before and after me},
but in {\it every} direction\ldots"

\vspace*{.5cm}
\large
\noindent
{\bf 5. Conclusion}

\vspace*{.1cm}
\noindent
We have outlined a conceptually very important extension of our pencil concept
of  the  time dimension in the case of Galois fields of characteristic two and
order greater than four. It has been  shown  that such a generalization of the
model may  not be a mere academic issue. On the contrary,  it seems to possess
a  serious  `observational/experimental'  counterpart  in  the  domain  of the
psychopathology  of  time.  The  issue obviously asks for and deserves further
effort and ingenuity to be properly explored and examined.

\vspace*{.5cm}
\noindent
{\bf Acknowledgement}

\vspace*{.1cm}
\noindent
The  work was partially supported by the 2001--2003 joint research project of
the Italian Research Council and the Slovak Academy of Sciences entitled `The
Subjective Time and its Underlying Mathematical Structure.'

\end{document}